\documentclass[12pt,thmsa]{article}%
\usepackage{amsfonts}
\usepackage{sw20bams}
\usepackage{amsmath}
\usepackage{amssymb}
\usepackage{graphicx}%
\providecommand{\U}[1]{\protect\rule{.1in}{.1in}}
\begin{document}

\author{Steven Finch}
\title{Idempotents and Nilpotents Modulo $n$}
\date{November 18, 2017}
\maketitle

\begin{abstract}
We study asymptotic properties of periods and transient phases associated with
modular power sequences. The latter are simple; the former are vaguely related
to the reciprocal sum of square-free integer kernels.

\end{abstract}

\footnotetext{Copyright \copyright \ 2006 by Steven R. Finch. All rights
reserved.}Let $\mathbb{Z}_{n}$ denote the ring of integers modulo $n$. Define
$S(x)$ to be the sequence $\{x^{k}\}_{k=0}^{\infty}$ for each $x\in
\mathbb{Z}_{n}$. We wish to understand the periodicity properties of $S(x)$,
that is, the statistics of
\[
\sigma(x)=
\begin{array}
[c]{l}%
\text{the \textbf{period}}\\
\text{of }S(x)
\end{array}
=
\begin{array}
[c]{l}%
\text{the least }m\geq1\text{ for which }x^{k+m}=x^{k}\text{ }\\
\text{for all sufficiently large }k,
\end{array}
\]
\[
\tau(x)=
\begin{array}
[c]{l}%
\text{the \textbf{transient}}\\
\text{\textbf{phase} of }S(x)
\end{array}
=
\begin{array}
[c]{l}%
\text{the least }\ell\geq0\text{ for which }x^{k+\sigma(x)}=x^{k}\text{ }\\
\text{for all }k\geq\ell.
\end{array}
\]
For example, the unique $x$ with $(\sigma,\tau)=(1,0)$ is $x=1$. If
$(\sigma,\tau)=(2,0)$, then $x$ is a square root of unity; if $(\sigma
,\tau)=(3,0)$, then $x$ is a cube root of unity \cite{FS}. If $\tau=0$ (with
no condition placed on $\sigma$), then $x$ is relatively prime to $n$. Hence
the number of such $x$ is
\[
\#\left\{  x\in\mathbb{Z}_{n}:x^{k}=1\text{ for some }k\geq1\right\}
=\varphi(n)
\]
where $\varphi$ is the Euler totient function and, asymptotically \cite{FS,
Apstl},
\[%
{\displaystyle\sum\limits_{n\leq N}}
\varphi(n)\sim\frac3{\pi^{2}}N^{2}=(0.303963550927...)N^{2}
\]
as $N\rightarrow\infty$. As another example, if $(\sigma,\tau)=(1,1)$, then
$x$ is an \textbf{idempotent}. The number of such $x$, including $0$ and $1$,
is
\[
\#\left\{  x\in\mathbb{Z}_{n}:x^{2}=x\right\}  =2^{\omega(n)}
\]
where $\omega(n)$ denotes the number of distinct prime factors of $n$ and
\cite{FS, Sec1}
\[%
{\displaystyle\sum\limits_{n\leq N}}
2^{\omega(n)}\sim\frac6{\pi^{2}}N\cdot\ln N
\]
as $N\rightarrow\infty$. More difficult examples appear in the following
sections. As in \cite{FS}, we make no claim of originality: Our purpose is
only to gather relevant formulas in one place.

\section{Generalized Idempotents}

\subsection{Bounded Transient Phase}

If $\tau\leq1$ (with no condition placed on $\sigma$), then the number of such
$x$ is \cite{Sec1}
\[
a(n)=\#\left\{  x\in\mathbb{Z}_{n}:x^{k+1}=x\text{ for some }k\geq1\right\}
.
\]
This is a multiplicative function of $n$ and
\[
a(p^{r})=\left\{
\begin{array}
[c]{lll}%
p &  & \text{if }r=1,\\
p^{r}-p^{r-1}+1 &  & \text{if }r\geq2.
\end{array}
\right.
\]
Let
\begin{align*}
F(s)  &  =%
{\displaystyle\sum\limits_{n=1}^{\infty}}
\frac{a(n)}{n^{s+1}}=%
{\displaystyle\prod\limits_{p}}
\left(  1+%
{\displaystyle\sum\limits_{r=1}^{\infty}}
\frac{p^{r}-p^{r-1}+1}{p^{r(s+1)}}\right) \\
\  &  =%
{\displaystyle\prod\limits_{p}}
\left(  1+\frac{p^{s+2}-2p+1}{p(p^{s+1}-1)(p^{s}-1)}\right)  =G(s)\cdot
\zeta(s).
\end{align*}
Hence, by the Selberg-Delange method \cite{FS, Slbrg, Dlnge, Tnbm1},
\[%
{\displaystyle\sum\limits_{n\leq N}}
a(n)\sim\frac12G(1)\cdot N^{2}=A\cdot N^{2}
\]
as $N\rightarrow\infty$, where
\[
A=\frac12%
{\displaystyle\prod\limits_{p}}
\left(  1-\frac1{p^{2}(p+1)}\right)  =0.440756919862...
\]
(the \textit{quadratic class constant} described in \cite{Finch}, divided by
two). Joshi \cite{Joshi} obtained this result via a different approach and
found an alternative formula:
\[
A=\frac3{\pi^{2}}%
{\displaystyle\sum\limits_{\ell=1}^{\infty}}
\left(  \frac1{\ell^{2}}\cdot%
{\displaystyle\prod\limits_{p|\ell}}
\frac p{p+1}\right)
\]
but did not numerically evaluate this expression.

If $\tau\leq2$ (with no condition placed on $\sigma$), then the number of such
$x$ is \cite{Sec1}
\[
b(n)=\#\left\{  x\in\mathbb{Z}_{n}:x^{k+2}=x^{2}\text{ for some }%
k\geq1\right\}  .
\]
This is a multiplicative function of $n$ and
\[
b(p^{r})=\left\{
\begin{array}
[c]{lll}%
p^{r} &  & \text{if }r\leq2,\\
p^{r}-p^{r-1}+p^{(r-1)/2} &  & \text{if }r\geq3\text{ and }r\equiv
1\operatorname*{mod}2,\\
p^{r}-p^{r-1}+p^{r/2} &  & \text{if }r\geq4\text{ and }r\equiv
0\operatorname*{mod}2.
\end{array}
\right.
\]
From
\begin{align*}%
{\displaystyle\sum\limits_{n=1}^{\infty}}
\frac{b(n)}{n^{s+1}}  &  =%
{\displaystyle\prod\limits_{p}}
\left(  1+%
{\displaystyle\sum\limits_{\substack{r\geq1, \\r\equiv1\operatorname*{mod}2
}}}
\frac{p^{r}-p^{r-1}+p^{(r-1)/2}}{p^{r(s+1)}}+%
{\displaystyle\sum\limits_{\substack{r\geq2, \\r\equiv0\operatorname*{mod}2
}}}
\frac{p^{r}-p^{r-1}+p^{r/2}}{p^{r(s+1)}}\right) \\
&  =%
{\displaystyle\prod\limits_{p}}
\left(  1+\frac{p^{3s+2}+p^{2s+2}-2p^{s+1}+p^{s}-2p+1}{p(p^{2s+1}%
-1)(p^{2s}-1)}\right)  ,
\end{align*}
we deduce that $\sum_{n\leq N}b(n)\sim B\cdot N^{2}$, where
\[
B=\frac12%
{\displaystyle\prod\limits_{p}}
\left(  1-\frac1{p^{2}(p^{2}+p+1)}\right)  =0.477176626987....
\]

If $\tau\leq3$ (with no condition placed on $\sigma$), then the number of such
$x$ is \cite{Sec1}
\[
c(n)=\#\left\{  x\in\mathbb{Z}_{n}:x^{k+3}=x^{3}\text{ for some }%
k\geq1\right\}  .
\]
This is a multiplicative function of $n$ and
\[
c(p^{r})=\left\{
\begin{array}
[c]{lll}%
p^{r} &  & \text{if }r\leq3,\\
p^{r}-p^{r-1}+p^{2(r-1)/3} &  & \text{if }r\geq4\text{ and }r\equiv
1\operatorname*{mod}3,\\
p^{r}-p^{r-1}+p^{(2r-1)/3} &  & \text{if }r\geq5\text{ and }r\equiv
2\operatorname*{mod}3,\\
p^{r}-p^{r-1}+p^{2r/3} &  & \text{if }r\geq6\text{ and }r\equiv
0\operatorname*{mod}3.
\end{array}
\right.
\]
From
\begin{align*}%
{\displaystyle\sum\limits_{n=1}^{\infty}}
\frac{c(n)}{n^{s+1}}  &  =%
{\displaystyle\prod\limits_{p}}
\left(  1+%
{\displaystyle\sum\limits_{\substack{r\geq1, \\r\equiv1\operatorname*{mod}3
}}}
\frac{p^{r}-p^{r-1}+p^{2(r-1)/3}}{p^{r(s+1)}}+%
{\displaystyle\sum\limits_{\substack{r\geq2, \\r\equiv2\operatorname*{mod}3
}}}
\frac{p^{r}-p^{r-1}+p^{(2r-1)/3}}{p^{r(s+1)}}\right. \\
&  \ \;\;\;\;\;\;\;\;\,\;+\left.
{\displaystyle\sum\limits_{\substack{r\geq3, \\r\equiv0\operatorname*{mod}3
}}}
\frac{p^{r}-p^{r-1}+p^{2r/3}}{p^{r(s+1)}}\right) \\
&  =%
{\displaystyle\prod\limits_{p}}
\left(  1+\frac{p^{5s+2}+p^{4s+2}+p^{3s+2}-2p^{2s+1}+p^{2s}-2p^{s+1}%
+p^{s}-2p+1}{p(p^{3s+1}-1)(p^{3s}-1)}\right)  ,
\end{align*}
we deduce that $\sum_{n\leq N}c(n)\sim C\cdot N^{2}$, where
\[
C=\frac12%
{\displaystyle\prod\limits_{p}}
\left(  1-\frac1{p^{2}(p^{3}+p^{2}+p+1)}\right)  =0.490145568004....
\]
The pattern exhibited by $A$, $B$, $C$ is clear and deserves proof for
$\tau\leq T$, for arbitrary $T$. A different attempt \cite{Patel1} to
determine the asymptotics of $\sum_{n\leq N}b(n)$ and of $\sum_{n\leq N}c(n)$
unfortunately turned out to be erroneous \cite{Patel2}.

\subsection{\label{idmpnts}Bounded Period}

If $\sigma=1$ (with no condition placed on $\tau$), then the number of such
$x$ is \cite{Sec2}
\[
u(n)=\#\left\{  x\in\mathbb{Z}_{n}:x^{k+1}=x^{k}\text{ for some }%
k\geq0\right\}  .
\]
This is a multiplicative function of $n$ and
\[
u(p^{r})=\left\{
\begin{array}
[c]{lll}%
2 &  & \text{if }r=1,\\
p^{r-1}+1 &  & \text{if }r\geq2.
\end{array}
\right.
\]
From
\begin{align*}%
{\displaystyle\sum\limits_{n=1}^{\infty}}
\frac{u(n)}{n^{s}}  &  =%
{\displaystyle\prod\limits_{p}}
\left(  1+%
{\displaystyle\sum\limits_{r=1}^{\infty}}
\frac{p^{r-1}+1}{p^{rs}}\right) \\
&  =%
{\displaystyle\prod\limits_{p}}
\left(  1+\frac{2p^{s}-p-1}{p(p^{s}-1)(p^{s-1}-1)}\right)
\end{align*}
we \textit{would like} to deduce that $\sum_{n\leq N}u(n)\sim U\cdot f(N)$ for
some simple expression $f(N)$. Unfortunately this is an unsolved problem. More
details are found in section [\ref{nilpnts}].

If $\sigma\leq2$ (with no condition placed on $\tau$), then the number of such
$x$ is \cite{Sec2}
\[
v(n)=\#\left\{  x\in\mathbb{Z}_{n}:x^{k+2}=x^{k}\text{ for some }%
k\geq0\right\}  .
\]
This is a multiplicative function of $n$ and
\[
v(p^{r})=\left\{
\begin{array}
[c]{lll}%
2^{r} &  & \text{if }p=2\text{ and }r\leq2,\\
2^{r-1}+4 &  & \text{if }p=2\text{ and }r\geq3,\\
3 &  & \text{if }p>2\text{ and }r=1,\\
p^{r-1}+2 &  & \text{if }p>2\text{ and }r\geq2.
\end{array}
\right.
\]
From
\begin{align*}%
{\displaystyle\sum\limits_{n=1}^{\infty}}
\frac{v(n)}{n^{s}}  &  =\left(  1+\frac2{2^{s}}+\frac4{2^{2s}}+%
{\displaystyle\sum\limits_{r=3}^{\infty}}
\frac{2^{r-1}+4}{2^{rs}}\right)  \cdot%
{\displaystyle\prod\limits_{p>2}}
\left(  1+%
{\displaystyle\sum\limits_{r=1}^{\infty}}
\frac{p^{r-1}+2}{p^{rs}}\right) \\
&  =\left(  1+\frac2{2^{s}}+\frac4{2^{2s}}+\frac{2^{s+2}-6}{2^{2s}%
(2^{s}-1)(2^{s-1}-1)}\right)  \cdot%
{\displaystyle\prod\limits_{p>2}}
\left(  1+\frac{3p^{s}-2p-1}{p(p^{s}-1)(p^{s-1}-1)}\right)
\end{align*}
we \textit{would like} to deduce that $\sum_{n\leq N}v(n)\sim V\cdot f(N)$,
where $f(N)\,$ is the same expression as for $\sigma=1$.

If $\sigma\leq3$ (with no condition placed on $\tau$), then the number of such
$x$ is \cite{Sec2}
\[
\#\left\{  x\in\mathbb{Z}_{n}:x^{k+3}=x^{k}\text{ or }x^{k+2}=x^{k}\text{ for
some }k\geq0\right\}  =w(n)-u(n)+v(n)
\]
where
\[
w(n)=\#\left\{  x\in\mathbb{Z}_{n}:x^{k+3}=x^{k}\text{ for some }%
k\geq0\right\}  .
\]
The latter is a multiplicative function of $n$ and
\[
w(p^{r})=\left\{
\begin{array}
[c]{lll}%
2 &  & \text{if }p=3\text{ and }r=1,\\
3^{r-1}+3 &  & \text{if }p=3\text{ and }r\geq2,\\
2 &  & \text{if }p\equiv2\operatorname*{mod}3\text{ and }r=1,\\
4 &  & \text{if }p\equiv1\operatorname*{mod}3\text{ and }r=1,\\
p^{r-1}+1 &  & \text{if }p\equiv2\operatorname*{mod}3\text{ and }r\geq2,\\
p^{r-1}+3 &  & \text{if }p\equiv1\operatorname*{mod}3\text{ and }r\geq2.
\end{array}
\right.
\]
From
\begin{align*}%
{\displaystyle\sum\limits_{n=1}^{\infty}}
\frac{w(n)}{n^{s}}  &  =\left(  1+\frac2{3^{s}}+%
{\displaystyle\sum\limits_{r=2}^{\infty}}
\frac{3^{r-1}+3}{3^{rs}}\right)  \cdot%
{\displaystyle\prod\limits_{p\equiv2\operatorname*{mod}3}}
\left(  1+%
{\displaystyle\sum\limits_{r=1}^{\infty}}
\frac{p^{r-1}+1}{p^{rs}}\right) \\
&  \cdot%
{\displaystyle\prod\limits_{p\equiv1\operatorname*{mod}3}}
\left(  1+%
{\displaystyle\sum\limits_{r=1}^{\infty}}
\frac{p^{r-1}+3}{p^{rs}}\right) \\
\  &  =\left(  1+\frac2{3^{s}}+\frac{2\cdot3^{s}-4}{3^{s}(3^{s}-1)(3^{s-1}%
-1)}\right)  \cdot%
{\displaystyle\prod\limits_{p\equiv2\operatorname*{mod}3}}
\left(  1+\frac{2p^{s}-p-1}{p(p^{s}-1)(p^{s-1}-1)}\right) \\
&  \ \cdot%
{\displaystyle\prod\limits_{p\equiv1\operatorname*{mod}3}}
\left(  1+\frac{4p^{s}-3p-1}{p(p^{s}-1)(p^{s-1}-1)}\right)
\end{align*}
we \textit{would like} to deduce that $\sum_{n\leq N}w(n)\sim W\cdot f(N)$,
where $f(N)\,$ is the same expression as for $\sigma=1$ and $\sigma\leq2$.

Note that $w-u+v$ is not multiplicative since $w(21)-u(21)+v(21)=8-4+9=13$
while $w(3)-u(3)+v(3)=2-2+3=3$ and $w(7)-u(7)+v(7)=4-2+3=5$. It would easily
follow that $\sum_{n\leq N}(w(n)-u(n)+v(n))\sim(W-U+V)\cdot f(N)$, completing
the case $\sigma\leq3$, \textit{if} the nature of $f(N)$ could be better ascertained.

\subsection{Unbounded Period}

Elements $x$ of small period are apparently quite rare for large $n$. We will
visit the other extreme. Consider, for example,
\[
m(n)=\#\left\{  x\in\mathbb{Z}_{n}:x^{k+\left\lceil \varphi(n)/2\right\rceil
}=x^{k}\text{ for some }k\geq1\right\}
\]
(the ceiling function is needed only for $1\leq n\leq2$, beyond which
$\varphi(n)$ is always even). This is not a multiplicative function, but
nevertheless can be simplified to
\[
m(n)=\left\{
\begin{array}
[c]{lll}%
3 &  & \text{if }n=4,\\
\dfrac{p^{r-1}(p+1)}2 &  & \text{if }n=p^{r},\text{ where }p>2\text{ and
}r\geq1,\\
p^{r-1}(p+1) &  & \text{if }n=2p^{r},\text{ where }p>2\text{ and }r\geq1,\\
n &  & \text{otherwise.}%
\end{array}
\right.
\]
From
\[%
{\displaystyle\sum\limits_{n=1}^{\infty}}
\frac{m(n)}{n^{s+1}}=\zeta(s)-\frac1{4^{s+1}}-\frac{2^{s}+1}{2^{s+1}}%
{\displaystyle\sum\limits_{p>2}}
\frac{p-1}{p(p^{s}-1)},
\]
we deduce that $\sum_{n\leq N}m(n)\sim(1/2)N^{2}$ since
\[
0<(s-1)%
{\displaystyle\sum\limits_{p}}
\frac{p-1}{p(p^{s}-1)}<(s-1)%
{\displaystyle\sum\limits_{p}}
\frac1{p^{s}}\sim-(s-1)\ln(s-1)\rightarrow0^{+}
\]
as $s\rightarrow1^{+}$. The behavior of $\sum_{n\leq N}(n-m(n))$ is more
subtle. From
\[%
{\displaystyle\sum\limits_{p}}
\frac{p-1}{p(p^{s}-1)}\sim-\ln(s-1)\sim%
{\displaystyle\sum\limits_{p}}
\frac1{p^{s}}
\]
and the fact that $\sum_{p\leq N}p\sim N^{2}/(2\ln(N))$ via the Prime Number
Theorem \cite{Odlyzk, BS}, we deduce that
\[%
{\displaystyle\sum\limits_{n\leq N}}
(n-m(n))\sim\frac38\frac{N^{2}}{\ln(N)}.
\]
It would be interesting to replace $\left\lceil \varphi(n)/2\right\rceil $ by
more slowly-growing expressions and to see what asymptotic consequences arise.

\section{\label{nilpnts}Nilpotents}

An element $x$ of $\mathbb{Z}_{n}$ is \textbf{nilpotent} if its power sequence
$S(x)$ is eventually zero. Define \cite{Sec3}
\[
z(n)=\#\left\{  x\in\mathbb{Z}_{n}:x^{k}=0\text{ for some }k\geq1\right\}  .
\]
This is a multiplicative function of $n$ and
\[
z(p^{r})=\left\{
\begin{array}
[c]{lll}%
1 &  & \text{if }r=1,\\
p^{r-1} &  & \text{if }r\geq2.
\end{array}
\right.
\]
Define also the \textbf{square-free kernel} $\kappa(n)$ to be the product of
all distinct prime factors of $n$. Clearly $\kappa(p^{r})=p$ for all $r\geq1$
and hence $z(n)=n/\kappa(n)$ for all $n\geq1$. From
\begin{align*}%
{\displaystyle\sum\limits_{n=1}^{\infty}}
\frac{z(n)}{n^{s}}  &  =%
{\displaystyle\prod\limits_{p}}
\left(  1+%
{\displaystyle\sum\limits_{r=1}^{\infty}}
\frac{p^{r-1}}{p^{rs}}\right) \\
\  &  =%
{\displaystyle\prod\limits_{p}}
\left(  1+\frac1{p(p^{s-1}-1)}\right)
\end{align*}
we \textit{would like} to deduce that $\sum_{n\leq N}z(n)\sim Z\cdot f(N)$ for
some simple expression $f(N)$. Unfortunately, as discussed in section
[\ref{idmpnts}], this is an unsolved problem. De Bruijn \cite{dBrjn1, dBrjn2,
dBrjn3, Tnbm2} proved that
\[
\ln\left(
{\displaystyle\sum\limits_{n\leq N}}
\frac1{\kappa(n)}\right)  \sim\left(  \frac{8\ln N}{\ln\ln N}\right)
^{1/2}\sim\ln\left(  \frac1N%
{\displaystyle\sum\limits_{n\leq N}}
\frac n{\kappa(n)}\right)
\]
and Schwarz \cite{Schwrz} proved that
\[%
{\displaystyle\sum\limits_{n\leq N}}
\frac1{\kappa(n)}\sim2^{-1/4}(4\pi)^{-1/2}\left(  \frac{\ln\ln N}{\ln
N}\right)  ^{1/4}\left(  \min_{0<y<\infty}N^{y}\cdot%
{\displaystyle\sum\limits_{n=1}^{\infty}}
\frac1{\kappa(n)n^{y}}\right)  .
\]
A\ more concrete rightmost factor would be good to see someday.

\section{Primitive Roots}

We have not mentioned the group $\mathbb{Z}_{n}^{*}$ (under multiplication) of
integers relatively prime to $n$ in this paper thus far. A well-known counting
problem concerns the number \cite{Sec4, Ng, IR, Ap}
\[
g(n)=\#\left\{  x\in\mathbb{Z}_{n}^{*}:\sigma(x)=\varphi(n)\right\}
\]
of \textbf{primitive }$\varphi(n)^{\text{th}}$\textbf{\ roots modulo} $n$.
Equivalently, $g(n)$ is the number of generators of $\mathbb{Z}_{n}^{*}$.
Clearly $g(n)>0$ if and only if $\mathbb{Z}_{n}^{*}$ is a cyclic group;
further,
\[
g(n)=\left\{
\begin{array}
[c]{lll}%
\varphi(\varphi(n)) &  & \text{if }n=1,2,4,\text{ }q^{j}\text{ or }%
2q^{j},\text{ where }q\text{ is an odd prime and }j\geq1,\\
0 &  & \text{otherwise.}%
\end{array}
\right.
\]
Also define the \textbf{reduced totient} or \textbf{Carmichael function}
\cite{EPS}
\[
\psi(n)=\left\{
\begin{array}
[c]{lll}%
\varphi(n) &  & \text{if }n=1,2,4\text{ or }q^{j},\text{ where }q\text{ is an
odd prime and }j\geq1,\\
\varphi(n)/2 &  & \text{if }n=2^{k},\text{ where }k\geq3,\\
\operatorname*{lcm}\left\{  \psi(p_{j}^{e_{j}}):1\leq j\leq\ell\right\}  &  &
\text{if }n=p_{1}^{e_{1}}p_{2}^{e_{2}}\cdots p_{\ell}^{e_{\ell}}\text{, where
}2\leq p_{1}<p_{2}<\ldots\text{ and }\ell\geq2,
\end{array}
\right.
\]
which is the size of the largest cyclic subgroup of $\mathbb{Z}_{n}^{*}$, and
consider the number \cite{Sec4, Li, CP, MSP}
\[
h(n)=\#\left\{  x\in\mathbb{Z}_{n}^{*}:\sigma(x)=\psi(n)\right\}
\]
of \textbf{primitive }$\psi(n)^{\text{th}}$\textbf{\ roots modulo} $n$. It is
known that \cite{Pillai}
\[%
{\displaystyle\sum\limits_{n\leq N}}
g(n)\sim\tilde A\frac{N^{2}}{\ln N}
\]
as $N\rightarrow\infty$, where
\[
\tilde A=\frac58%
{\displaystyle\prod\limits_{p}}
\left(  1-\frac1{p(p-1)}\right)  =0.233722383512...
\]
(five-eighths of \textit{Artin's constant} \cite{Finch}). A corresponding
result for $\sum_{n\leq N}h(n)$ evidently remains open. The issue of the
asymptotics of $\sum_{n\leq N}g(n)$ and of $\sum_{n\leq N}h(n)$ bears some
resemblance to the periodicity problems discussed earlier.

\section{Acknowledgements}

I thank G\'{e}rald Tenenbaum for suggesting the relation to $\sum_{n\leq
N}1/\kappa(n)$ and for informing me about \cite{dBrjn1, Schwrz}. I\ also thank
Pascal Sebah and Greg Martin, my coauthors in \cite{FS, FMS}. After this paper
was completed, I learned about \cite{LS}, which uses sophisticated tools to
examine mean periods over all $x\in\mathbb{Z}_{n}^{\ast}$ (a different and
more successful approach than mine). \ Subsequently, from \cite{Toth}, I
learned of recent progress \cite{RT} in understanding the asymptotics of $%
{\textstyle\sum\nolimits_{n\leq N}}
1/\kappa(n)$, which presumably carries over (in some manner) to both
$\sum_{n\leq N}u(n)$ and $\sum_{n\leq N}z(n)$. \ If someone should be so kind
as to work through the details, I\ would be grateful to hear about this.

\end{document}